\documentstyle[12pt]{article}
\newcommand{\no}{\noindent}
\newcommand{\nnb}{\nonumber}
\newcommand{\la}{\label}
\newcommand{\r}{r_{\varepsilon}}
\newcommand{\be}{\begin{equation}}
\newcommand{\e}{\end{equation}}
\newcommand{\ba}{\begin{eqnarray}}
\newcommand{\ea}{\end{eqnarray}}
\newcommand{\n}{\nabla}

\newcommand{\va}{\varepsilon}
\newcommand{\w}{\infty}
\newcommand{\1}{\forall t\in (0, \infty)}
\newcommand{\2}{\forall t\in [0, \infty)}
\newcommand{\f}{\frac}
\newcommand{\om}{\Omega}
\newcommand{\p}{\partial}

\begin{document}

\title {Translating Solitons  of Mean Curvature
 Flow of Noncompact Spacelike Hypersurfaces in Minkowski Space \thanks{Supported by
 National 973-Project from MOST and Trans-Century
Training Programme Foundation for the Talents by the Ministry of
Education. }}

\author{Huai-Yu Jian \\
  Department of  Mathematics\\
 Tsinghua University, Beijing 100084, P.R.China\\
(e-mail: hjian@math.tsinghua.edu.cn)}

\date {}

\maketitle

{\bf Abstract.}  In this paper, we study  the existence,
uniqueness  and  asymptotic behavior of rotationally symmetric
translating solitons of the mean curvature flow in Minkowski
space. We also study the  asymptotic behavior and the strict
convexity of
 general  solitons of such flows.



\section*{1. Introduction  and Main Results}
\setcounter{section}{1} \setcounter{equation}{0}

  Minkowski space $R^{n, 1}$ is the linear space $R^{n+1}$
endowed with the Lorentz metric
 $$ds^2=\sum_{i=1}^n dx^2_i - dx^ 2_{n+1}.$$
Spacelike hypersurfaces  in $R^{n,1}$ are Riemanian n-manifolds,
having an everywhere lightlike normal field $\nu$ which we assume
to be future directed and thus satisfy the condition $<\nu ,
\nu>=-1.$  Locally, such surfaces can be expressed as graphs of functions
$x_{n+1}=u(x_1, \cdots , x_n): R^n\longmapsto R$ satisfying the
spacelike conditions $|\n u(x)|<1$ for all $x\in R^n.$

If a family of spacelike embeddings $X_t=X(\cdot , t): R^n
\mapsto R^{n, 1}$ with corresponding hypersurfaces $M_t=X(R^n, t)$
satisfy the evolution equation
 \be \label{1.1}
\frac{\partial X}{\partial t}  =  H \nu   \e  on some time
interval, we say that the surfaces $\{M_t\}$ are evolved by {\sl
Mean Curvature Flow} (MCF). Here $H=div_{M_t} \nu$ denotes the
mean curvature  of the hypersurface $M_t.$  Let $V(\cdot , t)$ be
the graph expression of  $M_t.$  Then $|\n V(\cdot , t) |<1$ and
MCF equation (1.1) is equivalent, up to a diffeomorphism in $R^n,$
to the equation
 \be \label{1.2}
\frac{\partial V}{\partial t}  = \sqrt{1-|\nabla V|^2}
div(\frac{\nabla V}{\sqrt{1-|\nabla V|^2}}) \ \ in \ \ R^n.\e

MCF has been extensively studied in Euclidean space; see [1] and
the references therein, while in Minkowski space,  MCF was studied
in [2, 3] for compact hypersurfaces and in [4, 5] for noncompact
hypersurfaces.  The method of MCF was used in [2, 3] to
constructed spacelike hypersurfaces with prescribed mean
curvature, which, as it is well-known, have played  important
roles in studying Lorentzian manifolds. In particular, maximal
hypersurfaces, i.e., the ones with zero mean curvature, were used
by Schoen and Yau in the first proof of the famous positive mass
theorem [6].

The solutions of MCF (1.1) (or (1.2), equivalently)  which move by
vertical translation are called {\sl Translating Solitons}.
Therefore, a translating soliton of MCF (1.2) is characterized by
$V(x, t)=u(x)+t,$ where $u: R^n\mapsto  R$ is an initial spacelike
hypersurface satisfying \be \la {1.3} div(\frac{\nabla
u(x)}{\sqrt{1-|\nabla u(x)|^2}})= \frac{1}{\sqrt{1-|\nabla
u(x)|^2}}, \forall x\in R^n.\e The spacelike condition reads as
\be \la {1.4} |\n u(x)|<1, \forall x\in R^n.\e

Translating solitons can be regarded as a natural way of foliating
spacetimes by almost null like hypersurfaces. It may be expected
that this kind of translating solitons would have applications in
general relativity [5]. For this purpose, it is useful to
understand their geometric structure sufficiently. In [5], the
existence of smooth solutions of (1.3)-(1.4) was proved by a PDE
method. However, using  ODE techniques  we can  find  strictly
convex radially symmetric solutions of (1.3)-(1.4).

 \vskip 0.2cm {\bf Theorem 1.1} {\sl   There exists exact one solution $r\in C^2 [0, \infty )$ to initial
 value problem
 \be \la {1.5} \frac{r''(t)}{1-(r'(t))^2}
+\frac{n-1}{t}r'(t)=1, \ \ t\in (0, \infty ) \e and \be \label
{1.6} r(0)=r'(0) =0\e such that $u(x)=r(|x-x_0|)+u(x_0)$ in $R^n$
for any radially symmetric  $C^2$ solution $u$ of (1.3)-(1.4),
where $x_0$ is the vertex of $u.$ Moreover, the function $r\in
C^{\infty} [0, \infty )$ satisfies
 \be \la {1.7}
\frac{t}{\sqrt{n^2+t^2}}\leq r'(t)<1, \forall t\in [0, \infty )\e
and
\be \la {1.8}
 0< r''(t)\leq 1, \forall t\in [0, \infty ) . \e
 Therefore, all
 rotationally symmetric spacelike
translating solitons of MCF  (1.2)  is smooth, strictly
 convex, unique up to a translation in $R^{n+1},$ and of  linear growth.}

\vskip 0.2cm

To describe the asymptotic behavior of general  solitons as $|x|\to \w
,$  we use the tangent cones methods in [7, 8]  for entire
spacelike convex hypersurfaces  of constant mean curvature and in
[9] for constant Gauss curvature.
 Define the {\sl  blow down} of $F$
at infinity by
\begin{equation}\la {1.9}
V_F(x)=\lim_{\rho\to \infty}\f {F(\rho x)}{\rho}.
\end{equation}
Since
$\f {d}{d\rho}(\f {F(\rho x)}{\rho}-\f {F(0)}{\rho})\geq 0
 $ if $F$ is convex, and
$ \f {F(\rho x)}{\rho}-\f {F(0)}{\rho}\leq |x|$
 if $F$ is spacelike.    $V_F$ is
well-defined over $R^n$ and the limit in (1.9) is uniform  on any
compact set in $R^n $ if $F$ is a convex  function satisfying
(1.4). Using Theorem 1.1 and the methods in [7, 8], we will prove

 \vskip 0.4cm {\bf Theorem 1.2} {\sl Suppose that
  $ u$ is a  convex  solution to (1.3)-(1.4).
 Then  the blowdown function $V_u$
  is a positive homogeneous degree one convex function
satisfying the 1-Lipschitz condition
\begin{equation}\la {1.10}
|V_u(x)-V_u(y)|\leq |x-y|, \forall x, y \in R^n
\end{equation}
and the null condition, i.e., for any $x\in R^n$ and any $\delta
>0,$ there is
$y\in R^n$ such that
\begin{equation} \la {1.11}
|V_u(x)-V_u(y)|=|x-y|=\delta.
\end{equation}
\no Furthermore, one has
\begin{equation} \la {1.12}V_u(y)=\lim_{\rho\to \infty} \f {u (\rho y)}{\rho}
=1 \ \  uniformly \ \ for \ \  y\in \overline {\n u(R^n)}\bigcap
S^{n-1}
 \end{equation}and
\begin{equation}\la {1.13}V_u(x)=
\lim_{\rho\to \infty} \f{u (\rho x)}{\rho} =|x| \ \ for \ \  x\in
\overline {\n u (R^n)},
\end{equation}}
where  $\overline {\n u(R^n)}$ is the smallest closed set containing $\{ y: y=\n u(x), x\in R^n\}$
in $R^n.$

\vskip 0.3cm A natural question is whether any solution to
(1.3)-(1.4) is convex. This question seems very difficult to the
author. However, we obtain the following result which is related
to this question in some way.

\vskip 0.3cm {\bf Theorem 1.3}   {\sl  Let $u$ be a convex
solution of equation (1.3)-(1.4).  If  the set $\om _0 =\{x\in R^n
: (u_{ij}(x))>0\}$  is nonempty,  then $  \om _0=R^n .$}

\vskip 0.3cm  A similar result was obtained for the equation
$\Delta u=f(u, \nabla u)$ in [10], for the equation of entire
spacelike hypersurfaces of constant mean curvature in [8] and for
the mean curvature flow in Eclidean space in [11, 12,13].

  This paper is organized as follows. In
Section 2, we will use ODE theory and  a priori estimate
techniques to prove Theorem 1.1.  In Section 3, we will give the
proof of Theorem 1.2.    In Section 4, we will  prove Theorem 1.3.

\section*{2. Radially symmetric solutions }
\setcounter{section}{2} \setcounter{equation}{0}

We start with some simple facts which will be used throughout this section.

If  $u(x)=r(|x-x_0|)+u(x_0)$  and $u\in C^{k, \alpha }(R^n)$ for some $k\geq 1, 0\leq \alpha \leq 1$ with
$k+\alpha \geq 2, $  then $r\in  C^{k, \alpha }[0, \infty )$
since $r(t)=u((t, 0)+x_0)=u((-t,0)+x_0)$ for all $t\geq 0.$ Thus $r'(0)=0$ and equation (1.3) is equivalent to
 \be\label {2.1} \frac{r''(t)}{1-(r'(t))^2} +\frac{n-1}{t}r'(t)=1,
\forall t\in (0, \infty ) \e and
 \be \label {2.2} r(0)=r'(0) =0;\e
 the spacelike condition (1.4) is equivalent to
  \be \la {2.3} 0<r'(t)<1, \1 \e
and the strict convexity  to
 \be\label{2.4} 1\geq r''(t)>0, \2 .  \e

Conversely, if $r\in C^2 [0, \infty)$ is a solution to
(2.1)-(2.2), then  it follows from  a direct computation that $u(x)=r(|x|)\in C^{1,1}(R^n)$ is a
   solution to (1.3)-(1.4).  By the standard regularity theory of elliptic equations in [14] we
   see that $r(|x|)\in C^{\infty }(R^n)$  and thus $r\in C^{\infty }[0, \infty ).$

\vskip 0.2cm

{\bf Lemma 2.1}  {\sl   If  $r\in C^2 [0, \infty)$ is  a solution
to (2.1)-(2.4), then it satisfies
 (1.7).}

{\bf Proof.}
 If $r'(t)<1-\delta $ for all $t\in [0, \w )$ and some $\delta \in
 (0, 1),$ then  $r''(t)\geq \f {\delta }{2}$ for all $t\geq t_0$
 and for some large $t_0>0$ by (2.1). Integrating this inequality over
 $[t_0, t)$ we obtain
 $$ 1-\delta > r'(t) \geq \f {\delta }{2}(t-t_0)-r'(t_0)$$
 for all $t\geq t_0,$ a contradiction.  Therefore, there is a
 sequence $t_k\to \w $ such that $r'(t_k)\to 1.$  Using (2.4), we
 get
 \be \la {2.5} \lim _{t\to + \w } r'(t)=1.\e
 Note that the inequality on the right sides of (1.7) follows directly  from (2.3).  We want only to prove
\be \la {2.6} r'(t)\geq \f {t}{\sqrt{n^2+t^2}}, \forall t\geq 0.\e
 On the contrary that (2.6) is false. Then we have a $t_0>0$
such that
$$r'(t_0)< \f {t_0}{\sqrt{n^2+t_0^2}}.$$
Observing that $r'(0)=0$ and
$$\lim_{t\to +\w }(r'(t)-\f {t}{\sqrt{n^2+t^2}})=0$$
by (2.5), we see that the function $r'(t)-\f {t}{\sqrt{n^2+t^2}} $
attains its negative minimum at a point $t_1>0.$ Hence
$$ r''(t_1)=(\f {t_1}{\sqrt{n^2+t_1^2}})'=n^2(n^2+t_1^2)^{-\f {3}{2}}$$
and $$r'(t_1)<\f {t_1}{\sqrt{n^2+t_1^2}}.$$ This, together with
(2.1), imply
\[ \begin{array}{lll}
 1 &=&   \displaystyle\frac{r''(t_1)}{1-(r '(t_1))^2} +\frac{n-1}{t_1}r
'(t_1)\\
& < &  \displaystyle \f {n^2(n^2+t_1^2)^ {-\f {3}{2}}}{1-\f
{t_1^2}
{ n^2+t_1^2}}+\frac{n-1}{t_1}\cdot\f {t_1}{\sqrt{n^2+t_1^2}}\\
& = & \displaystyle \frac{n}{\sqrt{t_1^2+n^2}} <1,
\end{array}
\]
a contradiction!

 \vskip 0.4cm

{\bf Lemma 2.2} {\sl There exists a   $r\in C^{\infty }[0, \infty )$ to
(2.1)-(2.4).  }

{\bf Proof.}  Since  equation (2.1) is singular at  $t=0,$ we
consider the approximation problem
 \be \label {2.7}
\frac{r''(t)}{1-(r'(t))^2} +\frac{n-1}{t+\va }r'(t)=1, \forall
t\in (0, \infty )\e
 \be\label{2.8} |r'(t)|<1, \1 \e
 and
 \be \la {2.9} r(0)=0, \ \ r'(0)=\frac{\va }{n}
 .\e
Integrating (2.7) over $[0, t)$ we have
 $$\f {1}{2} [ \ln \f {1+r'(t)}{1-r'(t)}-\ln \f {1+r'(0)}{1-r'(0)}]
+(n-1)\int_0^t \f {r'(s)}{s+\va }ds =t,$$ which implies that for
any $R>0,$ there exist a constant $0<C(R)<1$ depending on $R$ such
that $$|r'(t)|<1-C(R), \forall t\in [0, R).$$ Therefore, by local
existence result of ODE, we see that for any $\va \in (0, 1)$
there is a unique smooth solution to (2.7)-(2.9). Denote this
solution by $\r .$  Obviously, \be \label{2.10} r_{\va}''(0)=
   [1-\f {n-1}{\va }r'(0)][1- (r'(0))^2]=\frac{n^2-\va ^2}{n^3}.\e
   This leads us to conclude
that \be \la {2.11} \r''(t)\geq 0, \ \
 \2 .\e
 Otherwise, there is a $t_1\in (0, \w )$ such that $\r
'' (t_1) < 0. $ Then we may choose $t_0>0$ and $\delta >0$ such
that
  \be \la {2.12} \r ''(t_0)=0, \r ''(t) <0 ,\forall t\in
(t_0, t_0+\delta ).\e By (2.9) and (2.10), we may further assume
\be \la {2.13} \r '(t)>0, \forall t\in [t_0, t_0+\delta ).\e
Hence, $ 0<\r '(t)<\r '(t_0)$ for all $t\in (t_0, t_0+\delta ).$
But this, together with (2.7), (2.8), (2.12) and (2.13), implies
$$1 =   \frac{n-1}{t_0+\va }\r '(t_0)
>  \frac{n-1}{t+\va }\r
'(t)  >  \frac{\r''(t)}{1-(\r '(t))^2} +\frac{n-1}{t+\va }\r '(t)
 =1$$ for all $t\in (t_0, t_0+\delta ),$
a contradiction! This proves (2.11).

It follows from (2.11) and (2.9) that
 \be \la {2.14}
\r'(t)\geq \frac{\va }{n},  \1 .\e Using this, (2.10) and (2.11)
again, we claim that \be \la {2.15} \r''(t)> 0, \ \
 \2 .\e
 In fact, on the contrary that there is a $t_2>0$
 such that $\r '' (t_2) = 0.$  Then the
  function
$$y(t):=\frac{n-1}{t+\va }\r
'(t)=1-\frac{\r''(t)}{1-(\r '(t))^2} $$ attains a maximum at
$t_2.$ Hence, $y'(t_2)=0$ and therefore,  $\r '(t_2)=0,$
contradicting with (2.14). This proves (2.15).

Now we use (2.8), (2.9), (2.14), (2.15)and (2.7) to see that
 \be
\la {2.16} \frac{\va }{n}\leq \r ' (t)<1, \2 \e
 \be
 \la {2.17} \frac{t\va  }{n}\leq \r (t)\leq t, \2\e
 \be
 \la {2.18}
0<\r ''(t) =  (1-\frac{n-1}{t+\va }\r '(t))(1-(\r '(t))^2)\leq
 1-\frac{(n-1)\va }{n(t+\va )}, \2 . \e
By  estimates (2.16)-(2.18) we can choose a subsequence $\va _k\to
0$ $(k\to \w )$  and a function $r\in C^{1, \alpha }[0, \w )$ (
$\alpha \in (0, 1)$ fixed) such that
 \be \la {2.19}r _{\va _k} \to
r \ \ in \ \
 C^{1, \alpha }[0, \w )\ \  as \ \ k\to \w .\e
 Obviously,
 \be \label{2.20} r(0)= 0=r'(0),\e
 \be \la {2.21} 0\leq r'(t)\leq 1\ \ and \ \ 0\leq r''(t)\leq 1, \2 .\e
 Furthermore, we can conclude that
  \be \la {2.22} 0\leq r'(t)<1, \2 \e
  Otherwise, there is a $t_3>0$ such that  $r'(t_3)=1$ and $0\leq
  r'(t)<1$
  for all $t\in [0, t_3).$ Integrating (2.7) for $r_{\va _k}$
  over $[\f {t_3}{2},
  t)$ we have
$$\f {1}{2} [ \ln \f {1+r_{\va _k}'(t)}{1-r_{\va _k}'(t) }-
\ln \f  {1+r_{\va _k}'(\f {t_3}{2})}{1-r_{\va _k}'(\f {t_3}{2})} ]
+(n-1)\int_{\f {t_3}{2}}^t \f {r_{\va _k}'(s)}{s+\va _k}ds =\f
{t_3}{2}, \forall t\in (\f {t_3}{2}, t_3).$$ Letting $k\to \infty
$ and $t\to t_3^-$ then, we obtain
$$ +\w -
\ln \f  {1+r '(\f {t_3}{2})}{1-r '(\f {t_3}{2})}  +2(n-1)\int_{\f
{t_3}{2}}^{t_3} \f {r '(s)}{s }ds =\f {t_3}{2},  $$ a
contradiction! This shows (2.22). Observing that $\r $ satisfies
equation (2.7), we use (2.19)-(2.22) to see that
 $r\in C^2[0, \w )$ satisfies (2.1) and (2.2), which implies  $r\in C^{\w }[0, \w )$
 as we have said in the beginning of this section.

 Therefore,
 in order to finish the proof of Lemma 2.2, we want only to prove
 \be \la {2.23} r'(t)>0, \1 \e and
 \be \la {2.24} r''(t)>0, \2 .\e
 In fact,  by (2.10) we have
 $r'' (0) =\frac{1}{n} .$  Then (2.23) follows from the fact that
  $ r'(0)=0$ and
 $ r''(t)\geq 0$ for all $t>0$ as in (2.20) and (2.21).

 If there were a $t_4\in (0, \w )$ such that
 $r''(t_4)=0,$ then  it follows from (2.21) that the function
 $$ Z(t):=\frac{n-1}{t}r'(t)=1-\frac{r''}{1-(r')^2}$$
 attains a maximum at the point $t_4.$
 Thus
 $$Z'(t_4)=0 \ \ and \ \ therefore \ \
 r'(t_4)=0,$$ contradicting
 (2.23).   This proves (2.24) and thus Lemma 2.2.
\vskip 0.2cm

{\bf Lemma 2.3}  {\sl If $r_1 , r_2 \in C^2[0, \w )$ are both
solutions to  initial problem (2.1), (2.2) and (2.3), then
$r_1(t)=r_2(t)$ for all $r\geq 0.$}

{\bf Proof.} Let $ u_i(x)=r_i(|x|) \ \ (i=1, 2).$  As we have seen, $u_i\in C^{\w }(R^n)$
are solutions of (1.3)-(1.4). Fix $t>0, $ arbitrarily.
 We see that both $u_1(x)$ and $u_2(x)+r_1(t)-r_2(t)$ are  solutions of the Dirichlet
 problem of equation (1.3)-(1.4) over the ball $B_t(0)$
 with the same boundary value $r_1(t).$  Thus   $u_1(x)=u_2(x)+r_1(t)-r_2(t)$ for all $x\in B_t(0)$
 by the uniqueness theorem  [14, Theorem 10.2]. Taking x=0, we obtain $r_1(t)=r_2(t).$

\vskip 0.2cm
 {\bf Proof of Theorem 1.1:}  observing the simple facts at the beginning of this section and
 using Lemmas 2.1, 2.2 and 2.3, we immediately   obtain Theorem 1.1.

 \section*{3.  Proof of Theorem 1.2}
\setcounter{section}{3} \setcounter{equation}{0}

In this section, we  use the concept of  tangent cones at infinity
to describe the asymptotic behavior of the solitons as $|x|\to \w
.$  This method was used in [7, 8] for entire spacelike convex
hypersurfaces  of constant mean curvature and in [9] for constant
Gauss curvature.

Recall that the
   $\underline {blow down}$ function  \be \la {3.1}
V_F(x)=\lim_{\rho\to 0\infty}\f {F(\rho x)}{\rho}\e
 is well defined over $R^n$ and    the
limit is uniform on any compact set in $R^n $ if $F$ is a convex
function satisfying (1.4).

  \vskip 0.4cm
  {\bf Lemma 3.1} {\sl
  If $u$ is a convex function satisfying (1.4),
then $V_u$ is a positively homogeneous degree one convex function
satisfying the 1-Lipschitz condition
\begin{equation}\la {3.2}
|V_u(x)-V_u(y)|\leq |x-y|, \forall x, y \in R^n;
\end{equation} while if  $u$ is a convex  solution to  (1.3)- (1.4),
then  $V_u$ satisfies the null condition, i.e., for any $x\in R^n$
and any $\delta
>0,$ there is
$y\in R^n$ such that
\begin{equation} \la {3.3}
|V_u(x)-V_u(y)|=|x-y|=\delta.
\end{equation}}

 {\bf Proof.}
 The convexity and the positive homogeneity
are obviously from the definition of $V_u$ and the convexity of
$u.$

For any $x, y \in R^n, $ by (1.4) we have
 $$
|V_u(x)-V_u(y)|\leq  \limsup_{\rho\to \infty} \f {|u(\rho
x)-u(\rho y)|}{\rho}   \leq  |x-y|.   $$ Hence, it is sufficient
to prove the null condition.   On the contrary that there would
exist an $x\in R^n$, $\delta >0$ and $\theta >0$ such that
$$V_u(y)\leq V_u(x)+(1-2\theta )\delta$$
for all $y\in R^n$ with $|x-y|=\delta.$ Observing that the limit
in (3.1) is uniform on any compact set, we may choose a $\rho_0>0$
so that \begin{equation} \la {3.4} u_\rho(y)\leq V_u(x)+(1-\theta
)\delta
\end{equation}for all $\rho>\rho_0$
and all $y\in B(x, \delta),$ where we have used the notation
$$B(x, \delta )=\{y\in R^n : |y-x|<\delta\} \ \  and \ \ u_\rho(x)=\f
{u(\rho x)}{\rho}.$$ It follows from (1.3)-(1.4) that $u_\rho$
satisfies \be \la {3.5}\displaystyle(
\delta_{ij}+\frac{(u_\rho)_i(x)(u_\rho)_j(x)}{1-|\n
u_\rho(x)|^2})(u_\rho)_{ij}= \rho, \ \ \forall x\in R^n\e and \be
\la {3.6} |\n u_\rho(x)|<1,\ \ \forall x\in R^n.\e

Let $r(|x|)$ be the same solution to (1.3)-(1.4) as in Theorem 1.1, where $r$ is the unique solution of
  (1.5) and (1.6). Then the
  function
$$ W(y)=W(y; \rho) :=V_u(x) + (\delta - \f {r(\rho \delta )}{\rho}+
\f {r(\rho |y-x|)}{\rho})-\theta \delta $$
 also satisfies the same
(3.5)-(3.6) as $u_\rho$ for any $\rho>0$ and any $x\in R^n.$ Note
that $$ W(y)= V_u(x)+(1-\theta )\delta, \forall y\in \partial B(x,
\delta).$$ We use (3.4) and maximum principle on the domain $B(x,
\delta)$ to obtain
$$u_\rho(y)\leq W(y; \rho), \forall y\in B(x, \delta).$$
Letting $\rho \to \infty ,$    we have
\[ \begin{array}{lll}
 V_u (x) &\leq&  \displaystyle  V_u(x)+ (\delta -\theta \delta)  +\lim_{\rho
 \to \w }
\f {r(\rho |y-x|)}{\rho}-\lim_{\rho\to \w } \f {r(\rho \delta )}{\rho}\\
& =&   V_u(x)+ (\delta -\theta \delta)+(|y-x|-\delta) \\
& = & V_u(x)+ |y-x| -\theta \delta.
\end{array}
\]
Here, in order to determine the
limit,  we have used  the estimate  $$\sqrt{n^2+t^2}-n\leq r(t)\leq t,$$ which follows directly from (1.6) and (1.7)  in Theorem 1.1.   Taking $y=x$ yields
$$V_u(x)\leq V_u(x)-\theta \delta, $$
a contradiction. In this way, we have  shown the desired lemma.
 \vskip 0.4cm

 Recall that the tangential mapping of convex function $V_u$
at a point $x_0\in R^n $ is defined by
$$T_{V_u}(x_0)=\{\alpha \in R^n: V_u(x)\geq \alpha \cdot
(x-x_0)+V_u (x_0), \forall x\in R^n\}.$$ Obviously, it is a
closed, convex set and equals to $\n V_u (x_0)$ if $V_u$ is
differential at $x_0.$  The tangent cone of $u$ is defined by
 $$T_{V_u}(R^n)=\bigcup_{x\in R^n}T_{V_u}(x).$$

{\bf Lemma 3.2} {\sl If $u$ is a convex function satisfying (1.4),
then its tangent cone satisfies
$$ \overline {T_{V_u}(R^n)}=T_{V_u}(0)=\overline {\n u(R^n)} .$$}

 {\bf Proof.}

   To show $T_{V_u}(0)\subset \overline {\n u(R^n)},$ we choose
$\xi \in T_{V_u}(0).$ Since $V_u(0)=0, $ $V_u(y)\geq \xi \cdot y$
for all $y\in R^n.$ Given a $\delta >0.$ Observing that the limit
$$ V_u(y)=\lim_{\rho\to 0}\f {u(\rho y)}{\rho}=\lim_{\rho\to 0}\f
{u(\rho x)-u(0)}{\rho}$$ holds uniformly on any compact set in
$R^n,$ we see that

$$ \phi (y):=\f {u(\rho_\delta
y)-u(0)}{\rho_\delta }-\xi \cdot y +|y|^2\geq \f {\delta ^2}{2}$$
for all $ y\in \partial B(0, \delta )$ and some large
$\rho_\delta>1.$ But $\phi (0)=0,$ so $\phi$ attains its minimum
at a point $x_\delta \in B(0, \delta).$ Thus
$$ \n \phi (x_\delta )=\n u (\rho_\delta x_\delta)-\xi +2x_\delta =0.$$
Letting $\delta \to 0$ we get $\xi \in \overline {\n u(R^n)} .$
Therefore,$T_{V_u}(0)\subset \overline {\n u(R^n)}.$

To finish the proof, we follow the arguments in [7, p.793]. Let
$\xi\in T_{V_u}(R^n).$ Then there is an $x\in R^n $ such that
$$V_u(\rho y)\geq \xi \cdot (\rho y- x)+V_u (x), \forall y\in R^n,
\forall \rho>0.$$ Dividing this inequality by $\rho,$ using the
homogeneity of $V_u$ and then letting $\rho\to \infty ,$ we get
$$V_u(y)\geq \xi \cdot y, \forall y\in R^n.$$
This means $\xi \in T_{V_u}(0).$ Thus,  $\overline
{T_{V_u}(R^n)}=T_{V_u}(0)$ Since $T_{V_u}(0)$ is closed.

Now for any $x\in R^n,$ the convexity implies
$$ u(\rho y)\geq \n u(x) \cdot (\rho y- x)+u (x), \forall y\in R^n,
\forall \rho>0.$$ Dividing this  by $\rho,$  and   letting
$\rho\to \infty ,$ we see that
$$V_u(y)\geq \n u(x) \cdot y, \forall y\in R^n,$$
which implies $\n u(x)\in T_{V_u}(0).$ Since $x $ is arbitrary and
$T_{V_u}(0)$ is closed, we conclude that
$$\overline {\n u(R^n)}\subset T_{V_u}(0)= \overline {T_{V_u}(R^n)} .$$
This proves the  Lemma.

\vskip 0.4cm {\bf Proof of Theorem 1.2:} Since we have Lemma 3.1,
 it is enough to prove  (1.12) and (1.13).

 Choose $ y\in \overline {\n u
(R^n)}. $ By Lemma 3.2, $y\in  T_{V_u}(0).$  Because of
$V_u(0)=0,$ we have
$$ V_u(y) \geq y\cdot y =|y|^2.$$
  On the other hand, Lemma 3.1 yields
$$ V_u(y)\leq |V_u(y)-V_u(0)|\leq |y|.$$
Thus, we have
\begin{equation}\la { 4.5}|y|^2\leq V_u(y)\leq |y|, \forall y\in \overline
{\n u(R^n)}.
\end{equation}
Hence (1.12) follows.  Note that
$$ \f {u (\rho x)}{\rho}=|x| \f {u (\rho |x|\cdot\f{x}{|x|})}{\rho|x|}$$
and $\f {x}{|x|}\in S^{n-1}$ for $x\not= 0.$ This, together with
(1.12), yields (1.13).

\vskip 0.4cm

\section*{4.  Proof of Theorem 1.3 }
\setcounter{section}{4} \setcounter{equation}{0}

On the contrary that there exists a $x_1\in R^n \backslash\om _0.$
We will derive a contradiction. We may assume $\om _0$ is nonempty
and connected. (Otherwise, we replace it by one of its connected
components ).  Then there exists a short segment
 $l\subset \om _0$ such that $\bar {l}\cap \p \om _0=\{x_1\}$
 and $\va _1=dist (l, \p \om )>0.$ Take $x_2\in l$ such that $B_{\va }
 (x_2)\subset \om _0$ for some $\va \in (0, \va _1).$ Translating
 the ball  $B_{\va } (x_2)$ along the line $l$ we come to a point $\bar{x}$
 where the ball and $\p \om _0$
 are touched at the first time. It follows that
 \be \la {4.1}\bar{x}\in R^n \backslash\om _0, \ \ \       B_{\va } (x_0)\subset \om _0 \ \ \  and \ \ \ \overline{ B_{\va } (x_0)}
 \cap \p \om _0=\{\bar{x}\}\e
 for some $x_0\in \om _0.$ Moreover,  the minimum eigenvalue $\lambda (x)$ of
 the Hessian $(u_{ij}(x))$ satisfies $\lambda (\bar{x})=0.$
  By a coordinate translation and rotation we may
 arrange that
 \be \la {4.2} \bar{x}=0, \ \ u(0)=0, \ \ \n u(0)=0 \ \ and \ \ u_{11}(0)=\lambda (0)=0.\e
 Thus, the origin $0\in  \p B_{\va } (x_0)  $ and
 \be \la {4.3} (u_{ij}(x))>0\ \ in \ \  B_{\va } (x_0).\e
  Rewrite  equation (1.3)    as
 \be \la {4.4} \Delta u = 1+A(|\n u|^2)u_i u_j u_{i j}\ \ in \ \ R^n ,\e
 where $A(t)=\f {1}{t-1}$ is analytic for $t\in (-1, 1).$
 Differentiating (4.4) twice with respect to $\f {\p}{\p x_1},$  we have
 \begin{eqnarray}\la {4.5}
\Delta u_{11}&=& 4A''u_l u_{l 1}u_m u_{m1}u_i u_j u_{i j}+2A' u_{m 1} u_{m 1}u_{i} u_j u_{i j}\nnb \\
 & +& 2A' u_{m } u_{m 11}u_{i} u_j u_{i j}+ 8A' u_{m } u_{m 1}u_{i1} u_j u_{i j}\nnb \\
&+& 4A' u_{m } u_{m 1}u_{i} u_j u_{i j1}  +2A u_{i11} u_j u_{i j}\nnb \\
&+& 2A  u_{i1} u_{j1} u_{i j} +4A u_{i1} u_{j} u_{i j1}\nnb\\
&+& A u_{i} u_{j} u_{i j11} \ \  in \ \  R^n .
\end{eqnarray}
 Since $u$ is analytic in $R^n,$ we expand $u_{11}$ at $x=0$ as a power series to obtain
 $u_{11}(x)=P_k(x)+R(x)$ for all $x\in \overline {B_{\va }(x_0)} $
 (one can choose a smaller $\va $ in advance if necessary),
 where $ P_k(x)$ is the lowest order term, which, by (4.2) and (4.3),
 is a nonzero homogeneous polynomial of degree $k$, and $R(x)$ is the rest.
 The convexity of $u$ yields $k\geq 2.$
 It follows from (4.3) that  $u_{i i}u_{11}-(u_{i1})^2>0$ in
 $B_{\va }(x_0).$
 Summing over $i$  we have
 \be \la {4.6} \Delta u u_{11}>\sum_{j=1}^n u_{j1}^2\geq  u_{i1}^2\e
 for each $i=1, 2, \cdots , n.$

 We claim that each $u_{i 1}$ is of order at least $\f {k}{2}.$
  Otherwise, we expand $u_{i1}$ at $x=0$ as a power series
  so that the lowest order term $h(x)$ must be a a nonzero homogeneous polynomial.
 Choose
 $$a=(a_1, a_2, \cdots , a_n)\in B_{\va }(x_0)\backslash \{x\in B_{\va }(x_0): h(x)=0\}$$
 so that the  segment
 $$L=\{ ta:  t\in (0, 1)\} \subset
  B_{\va }(x_0).$$    Now restricting (4.6) on $L,$ multiplying
  the both sides by $t^{-k}$ and then letting $t\to 0^+,$  we see the limit of the left-hand side
  of (4.6) is a nonzero constant multiplied by $\Delta u(0)$ which equals to  $1$ by (4.4),
  but the limit of the right-hand side is positive infinite.
  This is a contradiction.

  Therefore,  each $u_{i 1}$ is of order at least $\f {k}{2}.$   Hence $u_{ij1},$
  $u_{11i}$ and $ u_{11ij}$ are of order at least $\f {k}{2}-1 ,$ $k-1$ and $k-2$ respectively.
  Also note that each $u_{i}$ is of order at least 1 by (4.2).
  With these facts one can check that the right-hand side of equation (4.5)
   is of order at least of $k ;$ while the left-hand side, $\Delta u_{11},$ is
   either of order $k-2,$ or $\Delta P_k=0$ for all $x\in  B_{\va }(x_0).$
   Since the first case is impossible by comparing  the orders of the two sides, we obtain that $P_k$
   is a harmonic polynomial in $  B_{\va }(x_0).$

  We claim that  $P_k\geq 0$ for all $x\in  B_{\va }(x_0).$    Otherwise,
  there exists $a=(a_1, a_2, \cdots , a_n)\in B_{\va }(x_0)$ such that $P_k(a)<0.$ Then
  $$\frac{u_{11}(ta)}{t^k}=P_k(a)+\frac{R(ta)}{t^k}, \ \ \forall t\in (0, 1),$$
  which implies $\lim_{t\to 0^+}\frac{u_{11}(ta)}{t^k}=P_k(a)<0$ contradicting the fact
  that $u_{11}>0$ in $B_{\va }(x_0) $ (see (4.3)).

  Now we use the strong maximum principle to see that $P_k>0$ for all $x\in  B_{\va }(x_0).$   But $P_k(0)=0,$
 and  it follows from Hopf's lemma that $\f {\p P_k}{\p \nu}(0)<0,$ where $\nu $ is the unit outward normal
 to the sphere $\p B_{\va }(x_0).$ This means that the degree of $P_k$ is only one, contradicting the
 fact $k\geq 2.$
This contradiction proves  the theorem.

\vskip 0.3cm

{\bf Acknowledgement.}  The first version of this paper was
completed while I was at Harvard University as a visiting scholar
from Sept.,2000 to July, 2001. I  would like to thank the
 Department  of Mathematics for the hospitality and Professor S. T.
Yau for his advice and encouragement.

\vskip 0.5cm
\section* { References}
\begin{enumerate}
\itemsep -2pt

\item[1]  G. Huisken, Local and global behavior of hypersurfaces
moving by mean curvature, Proceeding of Symposia in Pure
Mathematics, {\bf 54} (1993), 175-191.
\item[2] K. Ecker and G. Huisken, Parabolic methods for the
construction of spacelike slices of prescribed mean curvature in
cosmological spacetimes, Comm. Math. Phys., {\bf 135} (1991), 595-613.
\item[3] G. Huisken and S. T. Yau, Definition of center of mass
for isolated physical system and unique foliations by stable
spheres with constant curvature, Invent. Math., {\bf 124}  (1996), 281-311.
\item[4]K. Ecker, On mean curvature flow of spacelike
hypersurfaces in asymptotical flat spacetimes, J. Austral. Math.
Soc. Ser. A, {\bf 55} (1993), 41-59.
\item[5]K. Ecker, Interior estimates and longtime solutions for
mean curvature flow of noncompact spaceike hypersurfaces in
Minkowski space,  J. Differential  Geom., {\bf 45} (1997), 481-498.
\item[6]R. Schoen and S. T.Yau, Proof of the positive mass theorem
I, Comm. Math. Phys.,  {\bf 65}  (1979), 45-76.
\item[7] H. I. Choi and A. E. Treibergs,   Gauss maps of
spacelike constant mean curvature  hypersurfaces of Minkowski
space,  J.  Differential Geom., {\bf 32}  (1990), 775-817.
\item[8] A. E. Treibergs,  Entire spacelike
hypersurfaces of constant mean curvature in Minkowski space,
Invent. Math.,  {\bf 66} (1982), 39-56.
\item[9] B. Guan, H. Jian and R. Schoen,   Entire spacelike hypersurfaces
of constant Gauss curvature in Minkowski space ,  J. rene angew.
Math. , (to appear).
\item[10] N. Korevaar and J. Lewis, convex solutions to
       certain equations have constant rank Hessian, Arch Rational
       Mech. Anal., 97(1987), 19-32.

\item [11] H. Jian, Q. Liu and X. Chen, Convexity and symmetry of
translating solitons in mean curvature flow,  Chin. Ann. Math.,
26B(2005), 413-422.

\item [12] C. Gui and H. Jian, Solitons of mean curvature and
symmetry of solutions of fully nonlnear elliptic equations,
Preprint.

\item [13] X. Wang, Convexity solutions to the mean curvature
flow, arixv: math. DG/0404326, preprint.

\item [14] D. Gilbarg and N. S. Trudinger, Elliptic partial
differential equations of second order,  2nd Edition, Springer-Verlag, 1983.
 \end{enumerate}

\end{document}